\documentclass[conference]{IEEEtran}
\IEEEoverridecommandlockouts
\usepackage{cite}
\usepackage{amsmath,amssymb,amsfonts}
\usepackage{algorithmic}
\usepackage{graphicx}
\usepackage{textcomp}
\usepackage{xcolor}
\usepackage{url}
\def\BibTeX{{\rm B\kern-.05em{\sc i\kern-.025em b}\kern-.08em
    T\kern-.1667em\lower.7ex\hbox{E}\kern-.125emX}}

\usepackage{amsthm}

\newtheorem{lemma}{Lemma}

\theoremstyle{remark}

\newtheorem{assum}{Assumption}

\usepackage[hidelinks]{hyperref}
\usepackage{eso-pic}

\newcommand{\ConfNotice}{%
A.~M.~Ali and L.~Tirel, ``A Receding Horizon Control for General Assembly
Line Balancing Problems,'' in \emph{Proc.\ 2026 IEEE International Conference
on Systems, Man, and Cybernetics (SMC)}, Bellevue, WA, USA, Oct.~2026.}

\AddToShipoutPictureFG*{%
  \AtPageLowerLeft{%
    \raisebox{34pt}{%
      \makebox[\paperwidth][c]{%
        \parbox{\textwidth}{%
          \color{red}%
          \rule{\textwidth}{0.4pt}\par
          \vspace{2pt}
          \raggedright
          \footnotesize
          \bfseries
          \ConfNotice
        }%
      }%
    }%
  }%
}

\begin{document}

\title{A Receding Horizon Control For General Assembly
Line Balancing Problems\\

\thanks{Ali M. Ali is with the Department of Mechanical and Aerospace Engineering, Carleton University, Ottawa, ON, Canada. Luca Tirel is an Avionics Systems Engineer, Leonardo S.p.A., Rome, Italy.
}}

\author{
\IEEEauthorblockN{Ali Mohamed Ali and Luca Tirel}
}

\maketitle

\begin{abstract}
This paper introduces a novel approach to the General Assembly Line
Balancing Problem (GALBP) by utilizing a receding horizon optimal
control framework. The proposed discrete model for the assembly line
offers a flexible representation, avoiding assumptions about specific
line configurations. The control actions aim to optimize the industrial
assembly line by minimizing the completion time while adhering to
constraints such as task precedence, workstation capacity, and resource
requirements. Control actions are represented through task assignment
and resource allocation matrices, assigning tasks to specific workstations
and assigning resources to workstations, respectively. The optimization
problem is formulated as a Mixed-Integer Nonlinear Programming (MINLP)
problem. The inherent robustness of the receding horizon approach
ensures optimal solutions for the assembly line, effectively adapting
to sudden changes. Numerical experiments demonstrate the robustness
and effectiveness of the proposed control synthesis in efficiently
distributing tasks and resources, minimizing the overall completion
time.
\end{abstract}


\section{Introduction}

\begin{figure}[ht]
\includegraphics[width=0.9\columnwidth]{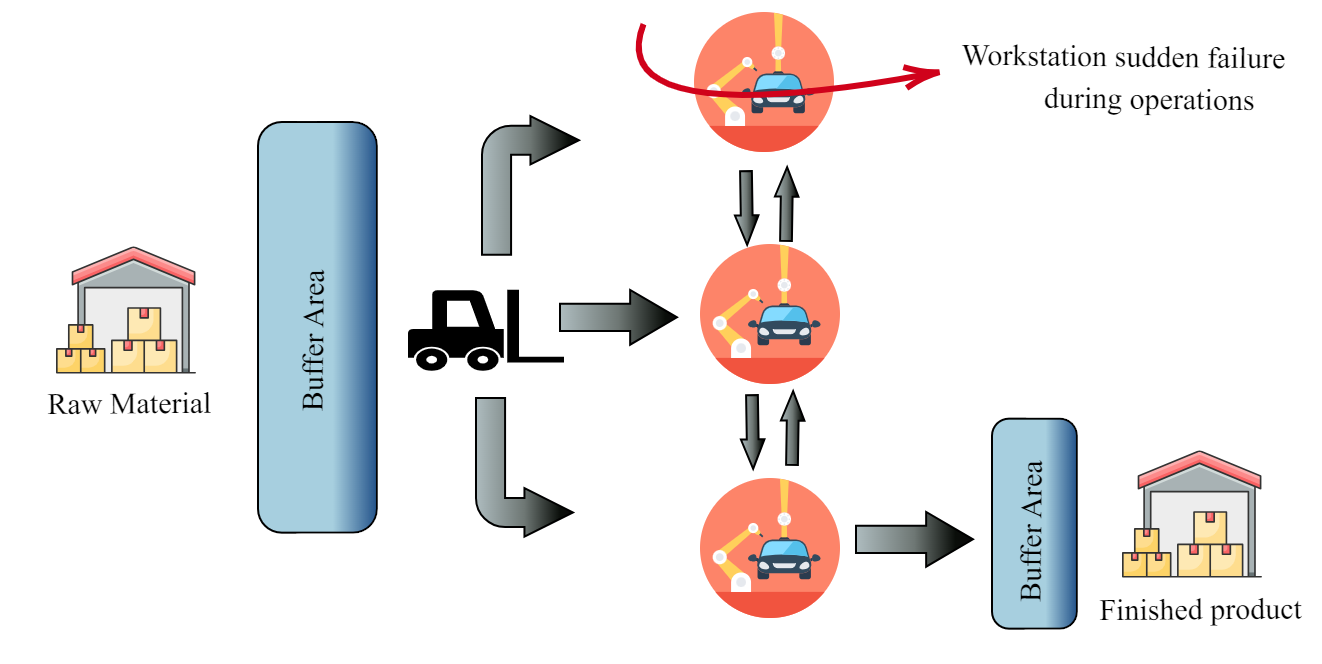}

\caption{A generic layout of an assembly line without any priori assumptions
of the shape of the line. Controlling the line involves assigning
tasks and needed resources to workstations to balance the workload,
and minimize the finishing time of the final product. The control
strategy should be robust against sudden changes in the line as a
failure of workstations.}
\label{fig:figure1}
\end{figure}

An assembly line is a manufacturing facility where individual components
of a product are assembled step-by-step to create a finished item.
This process plays a key role in allowing industries to efficiently
mass-produce goods \cite{sivasankaran2014literature}. The Assembly
Line Balancing Problem (ALBP), as described in the literature, involves
assigning tasks to a sequence of stations in a manner that satisfies
the precedence relationships between tasks while optimizing a performance
measure \cite{boysen2007classification}. Typically, obtaining an
exact solution for the ALBP is challenging because it requires solving
a complex nonlinear programming combinatorial optimization problem,
and being highly dependent on the specific configuration of the assembly
line. This complexity has been explored by a good deal of important
contributions \cite{salveson1955assembly,hazir2013assembly,tian2022enhanced,ali2023action}.
Unexpected changes in the assembly line (e.g., the sudden failure
of one or more workstations) may occur, requiring the problem to be
resolved and thereby provide optimal solutions. In this work, we propose
a receding horizon solution to the assembly line problem without any
a priori assumptions on the assembly line configuration (visit Fig.
\ref{fig:figure1}).

\paragraph*{Related Work}

The work in \cite{salveson1955assembly} presented an early attempt
for the ALBPs. Two main families of ALBP problems were identified
in \cite{salveson1955assembly}: Single ALBP (SALBP) and Generalized
ALBP (GALBP). The SALBP is characterized by single-element production,
serial line layout, a fixed common cycle time, a serial unilateral
line, precedence constraints, and deterministic operation times \cite{ali2025novel}.
The GALBP formulation relaxes some of the SALBP assumptions to allow
for more realistic scenarios. In particular, GALBPs allow for the
consideration of multi-model processes, zone constraints, delays,
parallel stations, and more complex layouts, to name a few. This latter
class of problems includes the U-shaped ALBP (UALBP) and the Mixed-Model
ALBP (MMALBP) \cite{kucukkoc2015balancing}. Another classification
of the ALBP is based on the objective function \cite{battaia2013taxonomy}.
Type-1 mainly focuses on minimizing the number of workstations
required for a fixed finishing time, and optimizing resource usage
while maintaining production rates. In contrast, Type-2 aims to minimize the finishing time for a given number of
workstations. Exact solutions
of SALBP are mainly based on repeatedly solving instances of the closely
related SALBP-1. In \cite{klein1996maximizing}, a branch and bound
algorithm is introduced. A multi-criteria decision-making approach
formulated and applied a goal programming model for U-type lines in
\cite{gokccen2006goal}. Motivated by the fact that exact solutions
are usually hard to compute, numerous research efforts have been directed
toward the development of computer-efficient approximation algorithms
or heuristics \cite{boysen2022assembly}. Numerous intelligent algorithms
have been extensively studied and applied providing sub-optimal solutions,
including the tabu search method \cite{lapierre2006balancing}, simulated
annealing method \cite{baykasoglu2006multi}, immune algorithm \cite{khoo2003line},
ant colony algorithm \cite{mcmullen2006multi,vilarinho2006antbal},
and genetic algorithm (GA) \cite{guo2008genetic,pistolesi2017emoga}.
Uncertainty and sudden changes may occur in the line, triggering recent
contributions toward robust solutions \cite{hazir2013assembly}. To
handle the uncertainty of activity times, a fuzzy-model-based solution
has been developed in \cite{ruppert2020fuzzy}. \textcolor{black}{The
aforementioned solutions are hand-crafted to specific types of uncertainty
and certain shapes of the line, for instance a U-shaped line with
uncertainty of activity times. This highlights the shortcomings of
the literature to compute exact solutions for the ALBP without priori
assumptions about the line layout and the nature of uncertainties
during operations.}

\textcolor{black}{\paragraph*{Contributions}}In this work, a receding
horizon approach is employed to solve ALBP, providing a robust solution
through its iterative implementation and effectively handling various
types of uncertainty, such as changes in available resources and workstation
failures.\textcolor{black}{{} The contributions of this work are
threefold. First, a discrete dynamical model for GALBP is designed
to fully capture the key features of the line while reducing the problem's
dimensionality to ensure computational feasibility. Second, a receding
horizon control approach is introduced to robustly solve the GALBP
against uncertainty and disturbances (e.g., sudden failure of one
or more workstations). Third, }the optimization problem within receding
horizon control is formulated as a Mixed Integer Nonlinear Programming
(MINLP) Problem. The computational complexity, feasibility, and closed-loop
stability of the proposed approach are investigated.

\textcolor{black}{\paragraph*{Structure}The remainder of the paper
is organized as follows: Section \ref{sec:Problem-Formulation} presents
mathematical notation, and problem formulation. Section \ref{sec:Sec5_Proposed-Scheme}
demonstrates the proposed control scheme. Section \ref{sec:Numerical-Results}
illustrates the effectiveness of the proposed scheme through numerical
simulations. Finally, Section \ref{sec:Sec7_Conclusion} concludes
the work.}

\section{\textcolor{black}{Problem Formulation \label{sec:Problem-Formulation}}}

\subsection{Notation}

In this paper, we refer to the sets of positive integers, binary numbers,
and the cardinality of a set as $\mathbb{Z}_{+}$, $\mathbb{Z}_{2}$,
and $|\cdot|,$ respectively. The time steps is denoted by $k\in\mathbb{Z}_{+}.$
The subscripts indicating workstations and tasks are $i$ and $j$,
respectively.

\subsection{Factory mathematical model}

A generic layout of a factory as shown in Fig. \ref{fig:figure1}
is described mathematically. The complexity of the description against
the dimensional scalability of the algorithms is a key issue that
arises as soon as the dimensions increase. This is the main reason
behind the explosion of computational time in the receding horizon
framework. The formal description of the factory is summarized in
the following items:

\paragraph{Tasks and Resources}

The set of all workstations is defined as $\mathcal{W}:=\{w_{1},\ldots,w_{I}\},$
while the set of all tasks are denoted as $\mathcal{T}:=\{T_{1},\ldots,T_{J}\}.$
We assume that there are different types of resources (e.g., raw materials,
electricity, workers, etc.) involved in the problem. The set of all
types of resources is $\mathcal{R}:=\{r_{1},\ldots,r_{R}\}$ with
the total number of types of resources being $R:=|\mathcal{R}|\in\mathbb{Z}_{+}$.
The total number of workstations and tasks are denoted by $I=|\mathcal{W}|\in\mathbb{Z}_{+}$
and $J=|\mathcal{T}|\in\mathbb{Z}_{+}$, respectively.

\paragraph{Constraints Description}

The main parameters describing the constraints are $O\in\mathbb{Z}_{+}^{I}$
and $U\in\mathbb{Z}_{+}^{R}$, which represent the workstation's
maximal capacities in terms of simultaneous processable tasks and
stored resources. Regarding the task constraints, we have $P\in\mathbb{Z}_{2}^{J\times J}$,
a Boolean matrix representing time precedence constraints between
the tasks (e.g., $T_{1}$ needs to be executed before $T_{2}$), $T_{j}\in\mathbb{Z}_{+}^{J}$
representing tasks' deadlines, $D\in\mathbb{Z}_{+}^{I\times J}$ that
indicates the duration of the tasks being dependent on the workstations
(e.g., $T_{1}$ can be executed in $w_{1}$ in two time steps while
$w_{2}$ can finish it in four time steps), and $C\in\mathbb{Z}^{J\times R}$
indicating the amount of specific resources needed for each task completion.
For example, a specific amount of $r_{1}$ and $r_{2}$ needs to be available
at a given time to execute $T_{1}$. Finally, there is $G\in\mathbb{Z}_{+}^{R}$
indicating the overall available resources for each type in the factory
(e.g., the available amount of raw material, electricity, workers,
etc.).

\paragraph{State Variables}

The number of tasks currently processed in the workstations also called
occupancies is indicated by $o[k]\in\mathbb{Z}_{+}^{I}$. For $I=3$,
$o[k]=[\begin{array}{ccc}
1 & 4 & 0\end{array}]$ indicates that $w_{2}$ is executing 4 tasks simultaneously at time
$k$. The remaining duration for the completion of the task in each
workstation is given by $d[k]\in\mathbb{Z}_{+}^{I\times J}$, and
the resources currently stored for the tasks' execution $r[k]\in\mathbb{Z}_{+}^{I\times J\times R}$.
Other state variables are needed to describe the start, execution,
and finishing of the tasks. $s[k]\in\mathbb{Z}_{2}^{I\times J}$ describes
whether a given task started at a specific workstation or not. Elements
of $s[k]$ remain equal to one if a specific task started at a workstation
and zero otherwise. Similarly, $e[k]\in\mathbb{Z}_{2}^{I\times J}$
and $f[k]\in\mathbb{Z}_{2}^{I\times J}$ describe the execution and
finishing of a specific task in a workstation, respectively. Elements
of $e[k]$ and $f[k]$ remain equal to one if a specific task started/finished
at a workstation and zero otherwise. Other auxiliary slack variables
are used to model the transition of the tasks between starting, execution,
and ending states as well as the availability of the resources in
the factory and workstations storage. $\text{Slack}_{1}[k]\in\mathbb{Z}_{2}^{I\times J}$
is needed for the dynamic duration update of a task performed by a
workstation at a given time. $\text{Slack}_{2}[k]\in\mathbb{Z}_{2}^{I\times J}$ describes the
time precedence constraints. $\text{Slack}_{3}[k]\in\mathbb{Z}^{I\times J}$
and $\text{Slack}_{4}[k]\in\mathbb{Z}^{I\times J}$ are used to update
the available resources in the factory and workstations' buffers with
the resources obtained from factory inventory, respectively. $\text{Slack}_{5}[k]\in\mathbb{Z}^{I\times J}$
is employed to remove resources from the workstations buffers when
a task ends. Finally, $\tau\in\mathbb{Z}_{+}$ is an auxiliary variable
of the objective function representing the finishing time of the last
task.

\paragraph{Control Variables}

Two control actions primarily concern the control variables. The
first is the assignment of a certain task to a specific workstation
while the second is the assignment of the resources from the main
factory storage inventory to the resources buffers located in different
workstations. $a_{i,j}^{T}[k]\in\mathbb{Z}_{2}^{I\times J}$ is the
control input for task assignment. Elements of $a_{i,j}^{T}[k]$ are
equal to one if a specific task is assigned at a workstation and zero
otherwise. For example, if $a_{i,j}^{T}[k](i,j)=1$, it implies that
the task $j$ is assigned to workstation $i$. Likewise, $a_{i,j,r}^{R}[k]\in\mathbb{Z}_{2}^{I\times J\times R}$
is the control input assignment for each type of the resources.

\subsection{Factory Dynamics }

In this subsection, the discrete dynamics describing how the state
and input variable evolve with time is introduced. The function $f(o[k],d[k],r[k],\ldots,\tau[k],a_{i,j}^{T}[k],a_{i,j,r}^{R}[k])$
that maps the current states at time $k$ to time $k+1$ is described.
Let the dynamics describing the occupancy level be expressed as follows:
\begin{equation}
o[k+1]=\begin{cases}
o[k](i)+1 & \text{if\hspace{3bp}}a[k](i,j)=1\\
o[k] & \text{otherwise}
\end{cases}\label{eq:1}
\end{equation}
When one task $j$ is assigned to workstation $i$, the corresponding
$i$-th entry of $o[k]$ is iteratively increased by one.\textcolor{black}{{}
Similarly, once a task is assigned, the remaining duration is equal
to the required time of the next task in the workstation. This information
is retained from the matrix $D$. As soon as the execution starts,
the correspondent duration decreases by one. Otherwise, it remains
constant. Next, the ($i$, $j$)-th entry of $d[k]$ is assigned to
the corresponding initial duration value taken from $D$. The dynamics
of the duration can be }described as follows:
\begin{equation}
d[k+1]=\begin{cases}
D(i,j) & \text{if\hspace{3bp}}a[k](i,j)=1\\
d[k](i,j)-1 & \text{if\hspace{3bp}}e[k](i,j)=1\\
d[k](i,j) & \text{otherwise}
\end{cases}\label{eq:2}
\end{equation}
\textcolor{black}{After the assignment, at each time step the duration
is decreased by one, transitioning tasks from execution state $e[k]\in\mathbb{Z}_{2}^{I\times J}$
to the finishing state $f[k]$. To minimize the dynamical constraints,
the resources action is kept Boolean, while the resources assigned
for task execution through the variable $r[k]$ are integer variables.
The resources assignment also corresponds to a reduction in the factory
inventories (e.g., $G[k+1](r)=G[k](r)-C(j,r)$). The resources assigned
are considered task-specific since this improves traceability and
allows modeling resources that do not vanish with consumption. The
resources dynamics can be }expressed as follows:
\begin{equation}
r[k+1]=\begin{cases}
r[k](i,j,r)+C(j,r) & \text{if}\ a_{i,j,r}^{R}[k]=1\\
r[k](i,j,r)-C(j,r) & \text{if}\ d[k](i,j)=0\\
r[k](i,j,r) & \text{otherwise}
\end{cases}\label{eq:3}
\end{equation}

\section{Proposed Receding Horizon Approach \label{sec:Sec5_Proposed-Scheme}}

This section details the proposed receding horizon control strategy.
The control objective is to optimally distribute the tasks and resources
among the workstations, reducing a given objective function. According
to receding horizon control, at each sampling time the current control
input is computed by solving a finite horizon open-loop optimal control
problem online, using the current state of the line (e.g., the current
occupancy level of each workstation $o[k]$) as the initial state of
the optimization problem. The first sample of the computed optimal
control sequence is applied to the plant, while the remaining part
of the control sequence is discarded. The open-loop optimal control
problem encodes the control objectives in the objective function and
includes dynamical and logical constraints such as task precedence.
In this work,
the focus will be on minimizing the finishing time of the final product
$\tau$. It is worth noting that the objective function
can be a weighted summation of terms that encodes other goals (e.g.,
minimization and balancing of the input/output resource inventory
in workstations). The optimal control problem that needs to be solved
for every time step starting from $k=0$ to $k=N$, where $N$ is
the prediction horizon, is presented as follows:{\small{}
\begin{equation}
\text{min}\hspace{0.5cm}\tau\label{eq:objective}
\end{equation}
subject to }{\small\par}

\hspace{3.5cm}\eqref{eq:1}-\eqref{eq:3}{\small{}
\begin{alignat}{1}
 & \sum_{k}\sum_{i}a_{i,j}^{T}[k]=1,\forall j.\label{eq:4}\\
 & d_{i,j}[k]a_{i,j}^{T}[k]=T_{ij}a_{i,j}^{T}[k],\forall k,\forall i,\forall j.\label{eq:5}\\
 & \tau\geq\sum_{k}a_{i,j}^{T}[k]*\left(k+T(i,j)\right),\forall k,\forall i,\forall j.\label{eq:6}\\
 & s_{i,j}[k]\geq\sum_{k}a_{i,j}[k],\forall k,\forall i,\forall j.\label{eq:7}\\
 & \sum_{i}s_{i,j}[k]\leq1,\forall k,\forall j.\label{eq:8}\\
 & f_{i,j}[k]\leq f_{i,j}[k+1],\forall k,\forall i,\forall j.\label{eq:9}\\
 & \sum_{k}d_{i,j}[k]=0,\forall i,\forall j.\label{eq:10}\\
 & a_{i,j}^{T}[k]e_{i,j}[k]=a_{i,j}^{T}[k],\forall k,\forall i,\forall j.\label{eq:11}\\
 & \sum_{i}e_{i,j}[k]\leq1,\forall k,\forall j.\label{eq:12}\\
 & f_{i,j}[k]+e_{i,j}[k]=s_{i,j}[k],\forall k,\forall i,\forall j.\label{eq:13}\\
 & e_{i,j}[k]\leq d_{i,j}[k],\forall k,\forall i,\forall j.\label{eq:14}\\
 & \text{Slack1}_{i,j}[k]=e_{i,j}[k],\forall k,\forall i,\forall j.\label{eq:15}\\
 & d_{i,j}[k+1]=d_{i,j}[k]-\text{Slack1}_{i,j}[k],\forall k,\forall i,\forall j.\label{eq:16}\\
 & d_{i,j}[k]\geq d_{i,j}[k+1],\forall k,\forall i,\forall j.\label{eq:17}\\
 & \text{Slack1}_{i,j}[k]=e_{i,j}[k-1]f_{i,j}[k],\forall k,\forall i,\forall j.\label{eq:18}\\
 & \sum_{k}\sum_{i}\text{Slack}2_{i,j}[k]=1,\forall j.\label{eq:19}\\
 & \text{Slack}2_{i,j_{1}}[k]=1\rightarrow\sum_{i}\sum_{k}a_{i,j_{2}}^{T}[k]=1,\forall k,\forall i.\label{eq:20}\\
 & \text{Slack}3_{r}[k]=\sum_{i}\sum_{j}a_{i,j,r}^{R}[k]T_{j,r}(j,r),\forall k,\forall i,\forall r.\label{eq:21}\\
 & \text{Slack}4_{r}[k]=\sum_{j}a_{i,j,r}^{R}[k]T_{j,r}(j,r),\forall k,\forall i,\forall r.\label{eq:22}\\
 & \text{Slack5}_{i,r}[k+1]=\sum_{j}e_{i,j}[k]f_{i,j}[k+1]T_{j,r}(j,r),\forall k,\forall i,\forall r.\label{eq:23}\\
 & a_{i,j,r}^{R}[k]=a_{i,j}^{T}[k],\forall k,\forall i,\forall j,\forall r.\label{eq:24}\\
 & o_{i}[k]=\sum_{j}e_{i,j}[k],\forall k,\forall i.\label{eq:25}\\
 & o_{i}[k]\leq W_{i}(i),\forall k,\forall i.\label{eq:26}\\
 & R[0]=R_{r}(r),\forall r.\label{eq:27}\\
 & R_{r}[k+1]=R_{r}[k]-\text{Slack3}[k+1],\forall k,\forall r.\label{eq:28}\\
 & r_{i,r}[1]=\text{slack4}{}_{i,r}[1],\forall i,\forall r.\label{eq:29}\\
 & r_{i,r}[k+1]=r_{i,r}[k]+\text{Slack4}_{i,r}[k+1]-\text{Slack5}_{i,r}[k+1],\forall k,\forall i,\forall r.\label{eq:30}
\end{alignat}
}{\small\par}

\paragraph*{Discussion on constraints}

The constraints within this model are structured to ensure proper
task assignment, resource allocation, and temporal management across
workstations. The constraint in \eqref{eq:4} ensures that each task
is assigned uniquely to a specific time instant on a workstation.
This is complemented by the constraint in \eqref{eq:5}, which involves
the duration assignment for each task, while the constraint in \eqref{eq:6}
serves as an auxiliary condition for the last time instant. Once the
starting state of a task is assigned, it is handled by the constraints
in \eqref{eq:7} and \eqref{eq:8}, ensuring that the state does not
reset. The constraints \eqref{eq:9} and \eqref{eq:10} impose conditions
on unique workstation assignment. Deadline enforcement is dictated
by constraint \eqref{eq:11}, whereas the constraints \eqref{eq:12}
and \eqref{eq:13} establish the execution mandate and the linkage
between the execution and finishing states. Constraint \eqref{eq:14}
ensures that execution only occurs when the duration is positive.
The model also defines slack variables through the constraints \eqref{eq:15},
\eqref{eq:18}, \eqref{eq:19}, \eqref{eq:21}, \eqref{eq:22}, and
\eqref{eq:23}, which govern slack definition, uniqueness, and time
precedence between tasks. The constraints in \eqref{eq:16} and \eqref{eq:17}
manage the update and the non-increasing nature of the task duration.
Furthermore, the constraints in \eqref{eq:24} to \eqref{eq:30} handle
the assignment of tasks and resources, the workstation occupancy levels,
and the resource initialization and dynamics, including the storage
of resources in the initial and subsequent time instants.

\begin{assum}\label{claim:Model assumpations}

\textit{The states of the assembly line are measurable which ensures
that the current occupancy level for all workstations are defined
at each time step.}

\end{assum}

Assumption \eqref{claim:Model assumpations} ensures that the feedback
signals sent to the receding horizon approach are measurable at each
time step. This guarantees that the receding horizon control will
re-optimize the line, even when new updates are introduced.
The optimization is a MINLP problem due to the presence of discrete
variables, such as task assignments and resource allocations, as presented
in constraints \eqref{eq:4}-\eqref{eq:13},\textcolor{black}{{} continuous
variables such as time management and resource dynamics in constraints
\eqref{eq:14} - \eqref{eq:30}, and nonlinear constraints as \eqref{eq:5}
and \eqref{eq:6}.}

\begin{figure}[h]
\begin{centering}
\centering\includegraphics[width=\columnwidth]{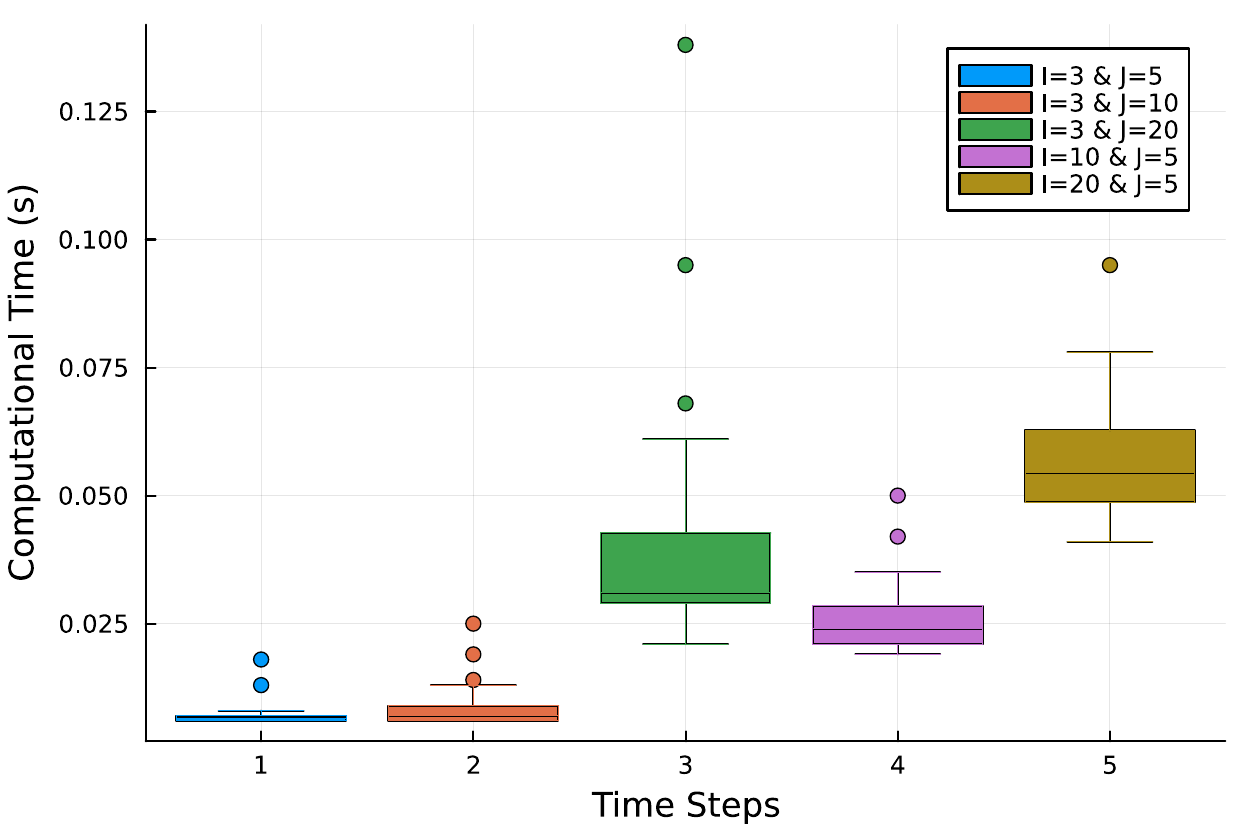}
\par\end{centering}
\caption{\label{fig:result_time} The average computational time for different
values of $I$ and $J$. }
\end{figure}

The receding horizon approach is implemented using the Julia programming
language \cite{bezanson2017julia} with the mixed-integer programming
solver Gurobi 11 \cite{gurobi}. The solver makes use of spatial
branch-and-bound and outer approximation \cite{vigerske2018scip}
as an alternative to the static piecewise-linear approximation of
the non-linear constraints. The computing resource used to carry out
the simulations in this work has the following specifications: Windows 11
(64 bit processor), CPU Intel(R) Core(TM) i7-10750H at 2.60 GHz with
6 cores and 12 logical processors, and 64 GB RAM. In Fig. \ref{fig:result_time},
the average computational time for different values of $I$ and $J$
is presented. One can notice the increase in computational time as
the number of workstations and tasks grows. It is worth noting that
the computational time is also a function of how tight the other
constraints are. For example, tighter constraints such as reducing
the deadline for completing tasks in $T_{j}$, or adding more time
precedence constraints between tasks in $T_{jj}$, will lead to an
increased computational time.

\paragraph*{Discussion on recursive feasibility}

The receding horizon iteration in \eqref{eq:objective} - \eqref{eq:30}
might be infeasible whenever there are conflicting constraints. For
instance, if the task deadline $T_{j}$ is less than the duration
that is needed by workstations to perform the task, the problem will
be infeasible. In this respect, the factory parameters need to be
carefully defined to avoid such infeasibility. In order to maintain
the recursive feasibility of the iterations at each time step, one
uses additional terminal constraints and terminal cost terms as indicated
by traditional receding horizon approaches \cite{rawlings2017model,fagiano2013generalized,ali2024mpc}.
However, introducing terminal constraints can be a challenging task
that also increases the computational complexity.

The proposed approach falls into the category of the economic receding
horizon approach. The following lemma provides the sufficient and
necessary condition of a general economic receding horizon approach
with an arbitrary economic objective.
\begin{lemma}
\label{lem:(Closed-Loop-stability}(Closed Loop stability of Economic
Receding Horizon \cite{diehl2010lyapunov}) If the Strong Duality
of steady-state solution and weak Controllability properties hold,
the steady-state solution of an arbitrary economic objective in the
closed-loop system will be asymptotically stable.
\end{lemma}

\paragraph*{Discussion on closed loop stability}

For more information about the formal definitions of the strong duality
and weak controllability properties please consult \cite{diehl2010lyapunov}.
In this work, closed-loop stability analysis is particularly challenging
due to the inclusion of Boolean variables. In addition, it differs
from the traditional receding horizon approaches which rely on a quadratic
objective function that can be treated as a Lyapunov function. In
the proposed approach the objective function cannot be directly used
as a Lyapunov function. This is a common difficulty associated with
the economic receding horizon approach. Moreover, the formulation
does not include specific terminal conditions to avoid limiting the
algorithm's feasibility region, as discussed previously. This omission
further complicates the stability analysis, where typical terminal
constraints are used to ensure stability under mild conditions \cite{rawlings2017model}.
Recalling Lemma \ref{lem:(Closed-Loop-stability}, an additional layer
of complexity arises when verifying the Strong Duality of the steady-state
solution to demonstrate closed-loop stability. For MINLPs, strong
duality does not always hold, and a duality gap may exist, where the
optimal values of the primal and dual problems differ. As such, the
closed-loop stability is deferred to future work.

\section{Numerical Results\label{sec:Numerical-Results}}

In this section, the proposed approach will address a simulation scenario
where one of the factory's workstations becomes unavailable for a
specific number of time steps. The full code for reproducing the
results is available at:
\nolinkurl{https://github.com/AliMohamedAliHassanAli/Receding-horizon-factory-Julia}.
The prediction horizon used in the simulation is $N=20$ while the
factory parameters are as follows: $I=3$, $J=5$, $R=2$, $O=\left[\begin{array}{ccc}
1 & 3 & 1\end{array}\right],$

\begin{flalign}
W_{\text{ir}} & =\left[\begin{array}{cc}
50 & 50\\
50 & 50\\
50 & 50
\end{array}\right],\label{eq:31}\\
T_{\text{ij}} & =\left[\begin{array}{ccccc}
4 & 6 & 8 & 5 & 2\\
5 & 6 & 12 & 3 & 6\\
3 & 4 & 5 & 3 & 5
\end{array}\right],\label{eq:32}\\
T_{\text{j}} & =\left[\begin{array}{ccccc}
20 & 20 & 20 & 20 & 20\end{array}\right],\label{eq:33}\\
T_{\text{jr}} & =\left[\begin{array}{cc}
13 & 12\\
6 & 7\\
4 & 3\\
2 & 4\\
3 & 5
\end{array}\right],\label{eq:34}\\
T_{\text{jj}} & =\left[\begin{array}{ccccc}
0 & 1 & 0 & 0 & 1\\
-1 & 0 & 0 & 0 & 0\\
0 & 0 & 0 & 0 & 0\\
0 & 0 & 0 & 0 & 1\\
-1 & 0 & 0 & -1 & 0
\end{array}\right],\label{eq:35}\\
R_{\text{r}} & =\left[\begin{array}{cc}
1000 & 2000\end{array}\right].\label{eq:36}
\end{flalign}

\begin{figure}[t]
\centering
\begin{minipage}{0.49\columnwidth}\centering
\includegraphics[width=\linewidth]{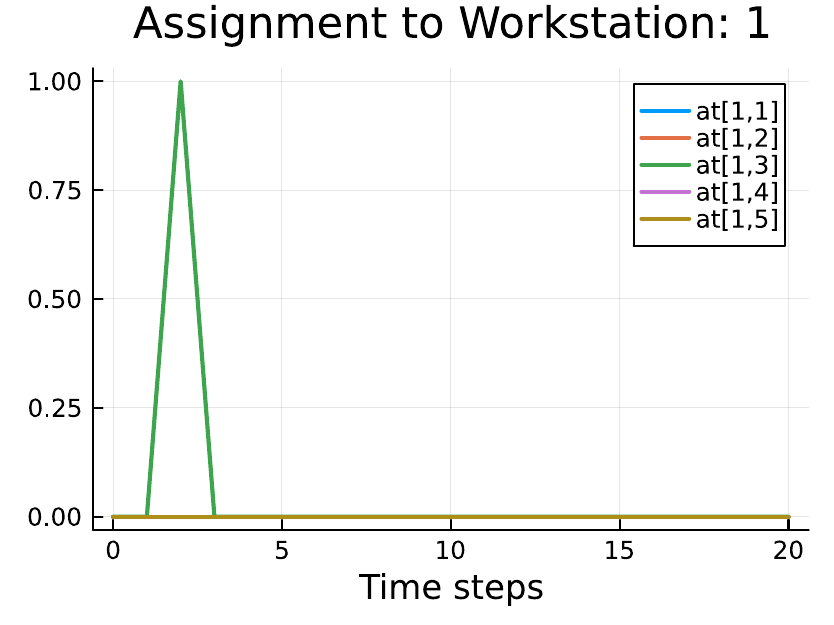}\\(a)
\end{minipage}\hfill
\begin{minipage}{0.49\columnwidth}\centering
\includegraphics[width=\linewidth]{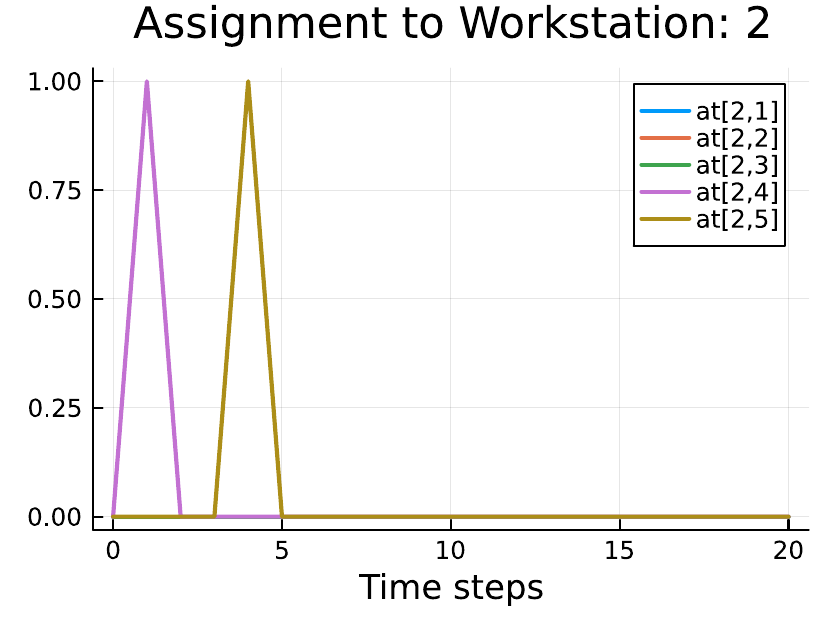}\\(b)
\end{minipage}

\vspace{4pt}
\begin{minipage}{0.62\columnwidth}\centering
\includegraphics[width=\linewidth]{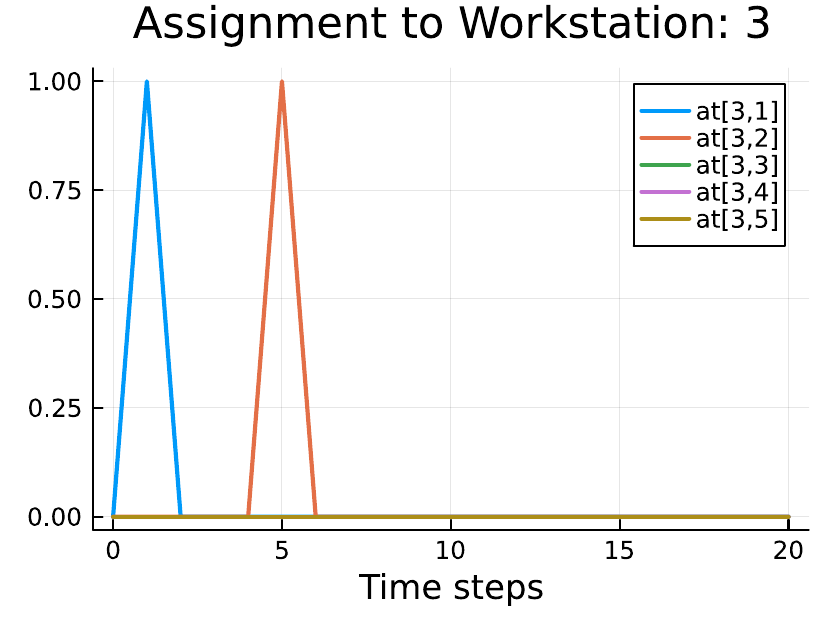}\\(c)
\end{minipage}
\caption{\label{fig:result1-1} Task assignment among the workstations. $T_{3}$
(Task number 3) is being assigned to $w_{1}$ (workstation 1) as shown
in item (a). While $T_{4}$ and $T_{5}$ are assigned to $w_{2}$ as
illustrated in item (b). Finally, the remaining two tasks $T_{1}$
and $T_{2}$ are assigned to $w_{3}$ as shown in item (c). }
\end{figure}
In the simulated scenario, $w_{3}$ will cease operation for 6
time steps starting from $k=8$. The proposed receding horizon approach
is expected to re-optimize the system to account for the disruption
caused by $w_{3}$ being out of service. Fig. \ref{fig:result1-1}
shows the task assignment to the workstations. One can verify that
the precedence constraints presented in \eqref{eq:35} are maintained.
The given $T_{jj}$ indicates that $T_{1}$ needs to be done before
$T_{2}$ and $T_{5}$ (this is clear from the first row where elements
in the second and fifth columns are 1 in \eqref{eq:35}). Moreover, from
\eqref{eq:32}, one can expect that the solver needs to assign $T_{1}$
to $w_{3}$ since the time needed to finish $T_{1}$ in $w_{3}$ is
3 while it is 4 and 5 in $w_{1}$ and $w_{2},$ respectively (see
the first column in \eqref{eq:32}). This is a logical assignment of
the tasks while minimizing the finishing time $\tau$.

\begin{figure}[t]
\centering
\begin{minipage}{0.49\columnwidth}\centering
\includegraphics[width=\linewidth]{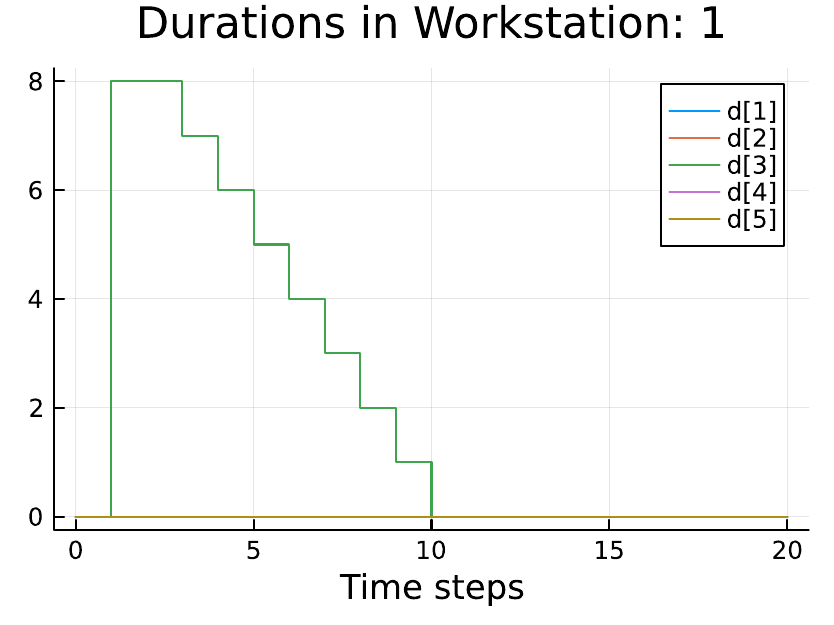}\\(a)
\end{minipage}\hfill
\begin{minipage}{0.49\columnwidth}\centering
\includegraphics[width=\linewidth]{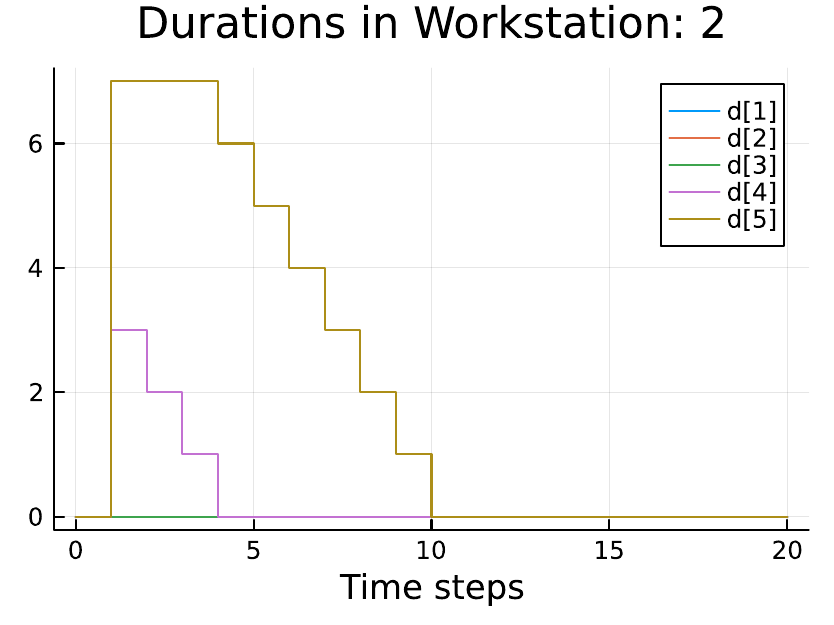}\\(b)
\end{minipage}

\vspace{4pt}
\begin{minipage}{0.62\columnwidth}\centering
\includegraphics[width=\linewidth]{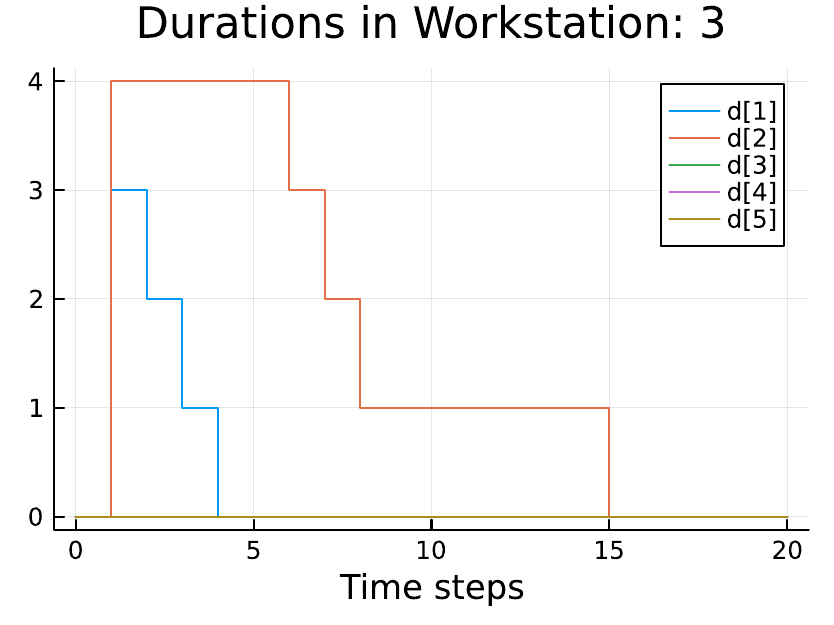}\\(c)
\end{minipage}
\caption{\label{fig:result2} Duration left for each Task in all workstations.
$T_{3}$ (Task number 3) in $w_{1}$ (workstation 1) finished at $k=10$
as shown in item (a). Item (b) illustrates the finishing time for
$T_{4}$ and $T_{5}$. Finally in item (c), one can notice the delay
in the duration left for $T_{2}$ in $w_{3}$ finishing at $k=15$
instead of $9$. }
\end{figure}

\begin{figure}[t]
\centering
\begin{minipage}{0.49\columnwidth}\centering
\includegraphics[width=\linewidth]{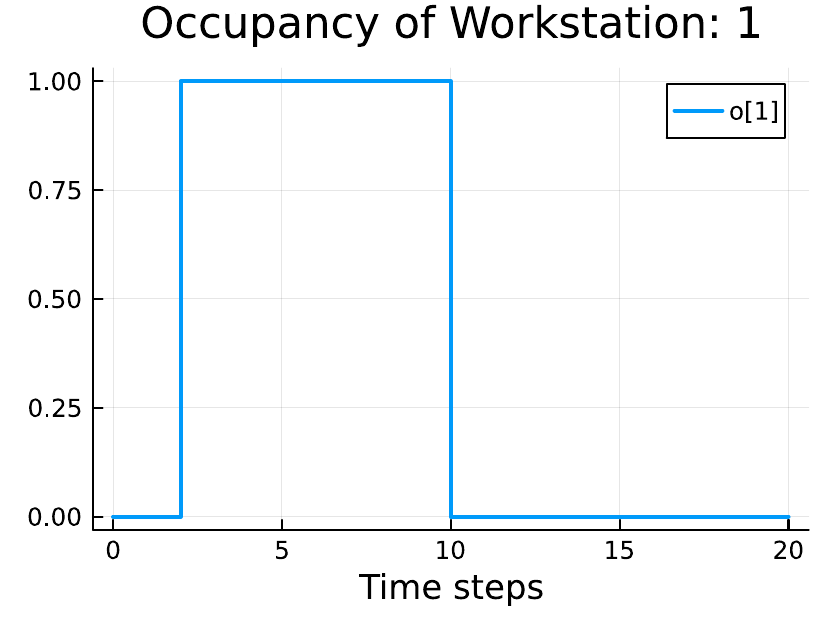}\\(a)
\end{minipage}\hfill
\begin{minipage}{0.49\columnwidth}\centering
\includegraphics[width=\linewidth]{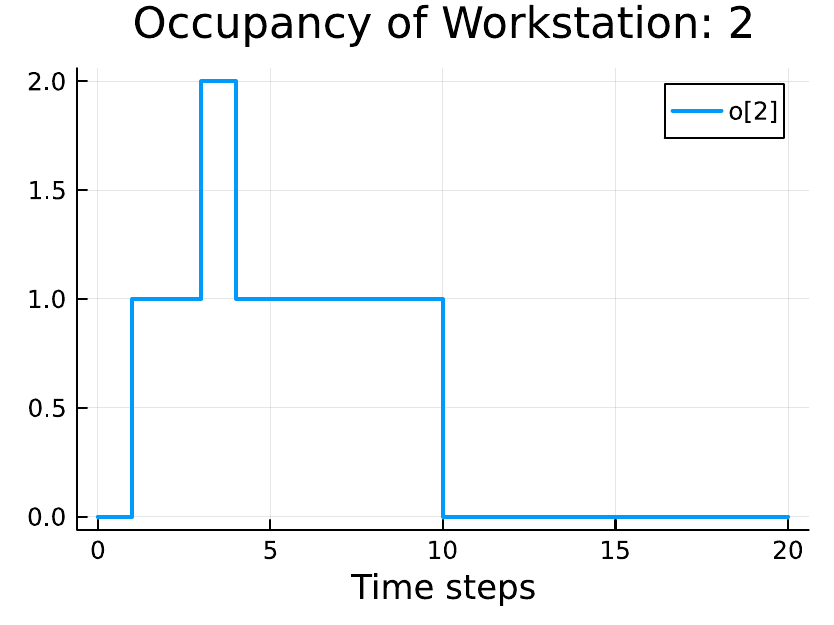}\\(b)
\end{minipage}

\vspace{4pt}
\begin{minipage}{0.62\columnwidth}\centering
\includegraphics[width=\linewidth]{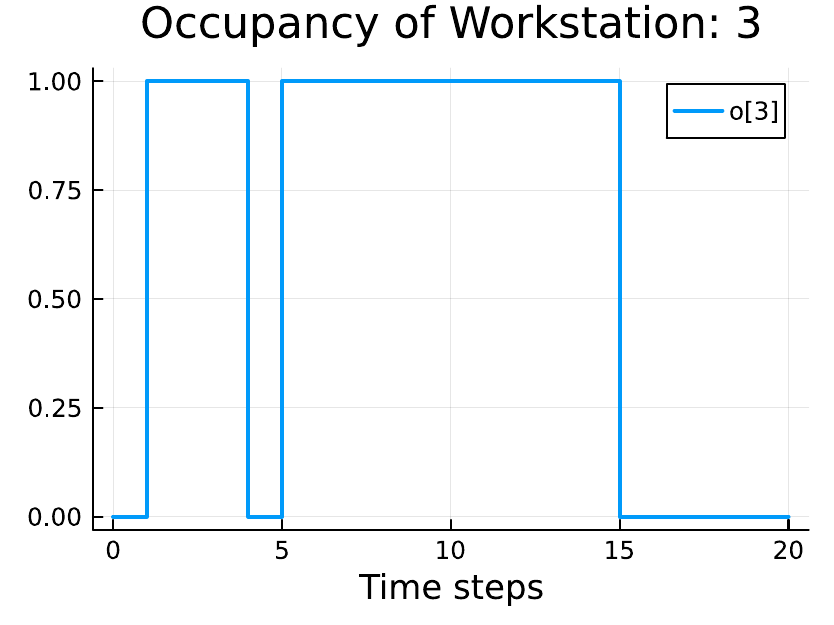}\\(c)
\end{minipage}
\caption{\label{fig:result3} Occupancy level $o[k]$ in each workstations.}
\end{figure}

Fig. \ref{fig:result2} shows the duration left for each task in every
workstation. As shown in Fig. \ref{fig:result1-1}, $T_{2}$ was assigned
to $w_{3}$ at $k=5$. According to \eqref{eq:32}, $T_{2}$ was expected
to be completed within 4 time steps. However, due to a 6-time-step
delay in $w_{3}$, the completion of $T_{2}$ was further postponed,
resulting in it finishing at $k=15$. In iterations after $k=10,$ the
receding horizon approach still assigned $T_{2}$ in $w_{3}$ despite
the delay in $w_{3}.$ To elaborate further, Fig. \ref{fig:result3}
depicts the occupancy level at each workstation. Similarly, the resources
assignment will make use of the initial resources inventory with values
in \eqref{eq:36}, to assign the needed resources for each type based
on \eqref{eq:34}. To conclude, the receding horizon approach assigned
optimally the tasks and resources resulting in a final finishing time
$\tau=15$ which is less than the deadlines of all the tasks (please
visit \eqref{eq:33}).

\section{Conclusion\label{sec:Sec7_Conclusion}}

\textcolor{black}{In this paper, we presented a receding horizon control
method to address the General Assembly Line Balancing Problem (GALBP).
By formulating the optimization challenge as a Mixed-Integer Nonlinear
Programming (MINLP) problem, we enhanced the system's robustness against
unexpected disruptions, such as workstation failures that could delay
product completion. Our approach showed strong performance in dynamically
adjusting task assignments and resource allocations, effectively minimizing
overall completion time while meeting operational constraints. The
closed-loop stability and recursive feasibility are partially investigated
and deferred to future work. Numerical experiments confirmed the effectiveness
of the proposed method, underscoring its potential to deliver robust
solutions in complex industrial environments. }

\vspace{0.25cm}

\bibliographystyle{IEEEtran}
\bibliography{ref}

\end{document}